\documentclass[12pt]{article}
\usepackage{geometry,color}                
\usepackage{float}
\usepackage{graphicx}
\usepackage{amssymb}
\usepackage{epstopdf}
\usepackage{amsmath}
\usepackage{marvosym}
\usepackage{mathtools}
\usepackage[round]{natbib}
\usepackage{multirow}
\usepackage{makecell}

\DeclareMathOperator*{\argmax}{arg\,max}

\newtheorem{theorem}{Theorem}
\newtheorem{conjecture}{Conjecture}

\newtheorem{proposition}{Proposition}

\newcommand{\Wbar}{{\overline{W}}}

\begin{document}
\def\spacingset#1{\renewcommand{\baselinestretch}%
	{#1}\small\normalsize} \spacingset{1}


{
	\title{\textbf{Cram\'er-von Mises tests for Change Points}}
	\author{Rasmus Erlemann\\ {\small NTNU, Department of Mathematical Sciences}\\
	Richard Lockhart\\
{\small SFU, Department of Statistics and Actuarial Science}\\
Rihan Yao\\
{\small SFU, Department of Statistics and Actuarial Science}}
	\maketitle
}

\thispagestyle{empty}

\begin{abstract}
	We study two nonparametric tests of the hypothesis that a sequence of independent observations is identically distributed against the alternative that at a single change point the distribution changes. The tests are based on the Cram\'er-von Mises two-sample test computed at every possible change point.  One test uses the largest such test statistic over all possible change points; the other averages over all possible change points. Large sample theory for the average statistic is shown to provide useful p-values much more quickly than bootstrapping, particularly in long sequences. Power is analyzed for contiguous alternatives.  The average statistic is shown to have limiting power larger than its level for such alternative sequences. Evidence is presented that this is not true for the maximal statistic.  Asymptotic methods and bootstrapping are used for constructing the test distribution. Performance of the tests is checked with a Monte Carlo power study for various alternative distributions. 
\end{abstract}

\noindent%
\vfill

\newpage

\section{Introduction}

Consider a sequence of independent observations $X_1,\ldots,X_n$.  We propose tests of the null hypothesis that the $X_i$ are independent and identically distributed (iid) with unknown continuous distribution $H$ against the change point alternative that there is 
some (unknown) $c$ with $ 1 \le c <n$ such that $X_1,\ldots,X_c$ are iid with continuous distribution $F$ and then $X_{c+1},\ldots,X_n$ are iid with some other continuous distribution $G$.  We will consider tests based on two sample empirical distribution function tests for equality of distribution, focusing on the two-sample Cram\'er-von~Mises test. 

If the time $c$ of the potential change point were specified in advance we could test the hypothesis that $F=G=H$ using any two sample test for equality of two distributions.  The two-sample Cram\'er-von~Mises test is one well known possibility.  Notation may be simpler to read if we used the shorthand $d=n-c$. Let
$$
F_c(x) =  \frac{1}{c} \sum_{i=1}^c 1(X_i \le x)
$$ 
be the empirical distribution function of the first $c$ observations and 
$$
G_d(x) =  \frac{1}{d} \sum_{i=c+1}^n 1(X_i \le x)
$$ 
be the empirical distribution function of the remaining $d$ observations. 
The combined empirical distribution function $H_{n}$ of the entire sample is 
$$
H_{n}(x) = \frac{cF_{c}(x)+d G_d(x)}{n}.
$$
The two-sample Cram\'er-von~Mises test of the hypothesis $F=G$ is based on the statistic
$$
W_{n}(c) = 
\frac{cd}{n} \int_{-\infty}^\infty \left\{F_c(x) -G_d(x)\right\}^2 dH_n(x).
$$
For a thorough discussion of this nonparametric test and a simple computing formula in terms of the ranks of the first $c$ values of $X$ in the whole sample see \cite{anderson1962}.  The distribution of the test statistic does not depend on $H$ under the null hypothesis provided $H$ is a continuous function.

A number of authors have suggested adapting this statistic to the change point problem. See, for instance, \cite{TimeSeriesChangePoint} and \cite{BookNonparametric} where the two natural possible test statistics considered herein are suggested and studied briefly. The first of these  tests can be used both to assess the existence of a change point and to estimate the location of the change if it exists. The statistic in question is
$$
W_{\max}\equiv \max_{1 \le c \le n-1} W_n(c).
$$
We shall also use $W_{\max}$ to define the estimated change point
$$\hat{c}_n = \argmax_{1\leq c\leq n-1}W_n(c) ;$$
thus $\hat{c}_n$ is the value of $c$ achieving the maximum. (We remark that the statistic $W_n$ is discrete and in small samples there is some modest probability that $\hat{c}_n$ will not be unique; this lack of uniqueness plays no role in the hypothesis testing problem.)

We prefer, however, the statistic
$$
\Wbar_n(X_1,\ldots,X_n) =\Wbar_n \equiv \frac{1}{n-1} \sum_{c=1}^{n-1} W_n(c).
$$
We offer several potential rationales for our choice:

\begin{itemize}
\item In many goodness-of-fit contexts quadratic statistics like ours outperform maximal statistics. For instance, the Cram\'er-von~Mises goodness-of-fit test is generally more powerful than the Kolmogorov-Smirnov test; see, for instance, \citet{stephensEDFCh4}.

\item Quadratic statistics such as we propose often have simpler large sample theory than do maximal statistics like the Kolmogorov-Smirnov test. Generally speaking the former have limiting distributions which are linear combination of chi-squares while the latter have limiting laws which are those of the supremum of a Gaussian process. The actual laws of these suprema are known only in special cases (although inequalities can often provide useful upper bounds on p-values).

\item The large sample theory in question often provides a more accurate approximation for quadratic statistics than it does for maximal statistics. For example, see \cite{KSlowpower} and \cite{Robustpower}.

\end{itemize}

In Section~\ref{sec:NullLimit} we present large sample distribution theory under the null hypothesis, show how to compute p-values based on this large sample theory and demonstrate that the asymptotic approximations are quite accurate for $n \ge 100$, particularly in the important lower tail.  Section~\ref{sec:PowerMC} presents a short power study showing that over a wide range of alternatives the statistic $\bar{W}$ is more powerful than $W_{\rm max}$.   Section~\ref{sec:ContiguousPower} presents asymptotic power calculations against contiguous sequences of alternatives; these permit useful approximations to the power of $\bar{W}$ in cases where  the null is not obviously false. By contrast, the limit theory for $W_{\rm max}$ does not lend itself to easy power calculations. We conjecture, however, that in this context of contiguous alternatives the statistic $W_{\rm max}$ has the defect that, unlike $\bar{W}$, its power converges to its level.  In this section we present some further Monte Carlo studies relevant to contiguous sequences of alternatives.
 Finally we present some discussion in Section~\ref{sec:Discussion}. We give proofs and evidence for the conjecture in the Appendix.

\section{Null limit theory}\label{sec:NullLimit}

Suppose that the null hypothesis holds and the $X_1,\ldots,X_n$ are iid with \emph{continuous} cdf $H$.  Then for all $c$ we have
$$
W(X_1,\ldots,X_{c},X_{c+1},\ldots,X_n)=W(H(X_1),\ldots,H(X_c);H(X_{c+1}),\ldots,H(X_n)).
$$

Thus in computing distribution theory under the null we may, and will, assume that $H$ is the uniform distribution; to emphasize the point we let $U_1,U_2,\cdots$ be an iid sequence of Uniform random variables; the joint law of $\left(H(X_1),\ldots,H(X_n)\right)$ is the same as that of $\left(U_1,\ldots,U_n\right)$.

Large sample theory for the two sample Cram\'er-von~Mises statistic is well known: 
if $c$ depends on $n$ in such a way that $c/n \to s \in (0,1)$ (or even just $\min\{c,n-c\}\to\infty$) then  
$$
W_n(c) \Rightarrow \sum_{j=1}^\infty \frac{Z_j^2}{\pi^2 j^2}
$$
where the $Z_i$ are iid standard normal; see \cite{anderson1962}. (Notice that the limit is free of $s$.) Our statistic has a related limit given as follows.

\begin{theorem} \label{theorem:wbar}
As $n \to \infty$ we have, under the null hypothesis,
$$
\Wbar_n \Rightarrow \Wbar_\infty \equiv \sum_{j=1}^\infty \sum_{k=1}^\infty \frac{Z_{jk}^2}{j(j+1) \pi^2 k^2}
$$
 where the  $Z_{jk}$ are iid standard normal. 
 \end{theorem}

The theorem is a consequence, as usual, of a suitable weak convergence result which we now present; the Gaussian process limit we derive is mentioned in \cite{TimeSeriesChangePoint}; the specific weights in Theorem~\ref{theorem:wbar} do not seem to have been previously described.  

We begin by defining the partial sum empirical process \cite[p. 225]{van1996weak}, 
for $(s,t) \in [0,1]^2$, by
$$
\mathbb{Z}_n(s,t) = \frac{1}{\sqrt{n}} \sum_{1 \le i \le ns} \left\{1(U_i \le t) - t\right\}.
$$
Our statistic can be described in terms of this process. 
Notice that 
$$
F_c(t) = \frac{\sqrt{n}}{c} \mathbb{Z}_n(c/n,t)  + t
$$
and that
$$
G_d(t) = \frac{\sqrt{n}}{d} \left\{\mathbb{Z}_n(1,t)-\mathbb{Z}_n(c/n,t)\right\} + t.
$$
Thus 
$$
F_c(t) -G_d(t) = \sqrt{n}\left\{\frac{ \mathbb{Z}_n(c/n,t)}{c} -\frac{\mathbb{Z}_n(1,t)-\mathbb{Z}_n(c/n,t)}{d}\right\}.
$$
We now define a process $\mathbb{W}_n(s,t)$ for $0<s<1$ and $0 \le t \le 1$ by
$$
\mathbb{W}_n(s,t) = \sqrt{s(1-s)} \left\{ \frac{\mathbb{Z}_n(s,t)}{s} - \frac{\mathbb{Z}_n(1,t)-\mathbb{Z}_n(s,t)}{1-s}\right\}
=\frac{\mathbb{Z}_n(s,t)-s\mathbb{Z}_n(1,t)}{\sqrt{s(1-s)}}.
$$
For given $c$ our  two sample test statistic is given by 
$$
W_n(c) = \int_0^1 \left\{\mathbb{W}_n(c/n,t)\right\}^2\, dH_n(t).
$$
The processes $\mathbb{Z}_n$ and $\mathbb{W}_n$ have well known weak limits given the in following theorem. It will also prove useful to introduce the notation
$$
\mathbb{B}_n(s,t) = \mathbb{Z}_n(s,t)-s\mathbb{Z}_n(1,t).
$$

\begin{theorem} \label{theorem:weakconvergence}
Under the null hypothesis: 
\begin{enumerate}
    
    \item As $n\to\infty$,
    $$
    \mathbb{Z}_n(s,t)\rightsquigarrow \mathbb{Z}_\infty
    $$
    a mean 0 Gaussian Process with covariance function 
    $$
    \rho_Z(s,t;s',t') = s\wedge s' \psi(t,t')
    $$
    where $\psi(t,t') = t\wedge t' - tt'$;
    
    \item As $n\to\infty$,
    $$
    \mathbb{B}_n(s,t)\rightsquigarrow \mathbb{B}_\infty
    $$
    a mean 0 Gaussian Process with covariance function
    $$
    \rho_B(s,t;s',t') = \psi(s,s') \psi(t,t');
    $$
    
    \item As $n\to\infty$,
    $$
    \mathbb{W}_n(s,t)\rightsquigarrow \mathbb{W}_\infty
    $$
    a mean 0 Gaussian Process with covariance function
    $$ 
    \rho_W(s,t;s',t') = \chi(s,s') \psi(t,t')
    $$
    where 
    $$
    \chi(s,s') =\frac{\psi(s,s')}{\sqrt{s(1-s)s'(1-s')}}.
    $$
    \end{enumerate}

\end{theorem}

The process $\mathbb{B}$ is called a Brownian pillow by some writers or a 4 side tied down Brownian motion; see, for instance \cite{zhang2014} or \cite{mckeague}. The process $\mathbb{Z}$ is a Blum-Kiefer-Rosenblatt process ; see \cite{blum1961}.  

We now record well known facts about the eigenvalues of the covariance $\rho_W$.
The covariance kernel $\psi$ is that of a Brownian Bridge. It has eigenvalues of the form $1/(\pi^2 k^2)$ for $k=1, 2,\cdots$ with corresponding orthonormal eigenfunctions $f_{\psi,k}(u) = \sqrt{2} \sin(\pi k u)$. The covariance kernel $\chi$ arises in the study of the Anderson-Darling goodness-of-fit test. It has eigenvalues of the form $1/\{j(j+1)\}$ 
 for $j=1, 2,\cdots$.  The corresponding orthonormal eigenfunctions are associated Legendre functions. The 
 $j^{\rm th}$
 eigenfunction is
 $$
 f_{\chi,j}(u) = 2\sqrt{\frac{2j+1}{j(j+1)}} \sqrt{s(1-s}) q_j(2s-1)
 $$
 where the $q_j$ are polynomials of degree $j-1$ defined recursively as follows:
 $$
 q_1(u) = 1,
 $$
 $$
 q_2(u)=3u
 $$
 and for $ j \ge 2$ 
 $$
 q_{j+1}(u) = \frac{1}{j} \left\{(2j+1) u q_j(u)-(j+1)q_{j-1}(u)\right\}.
 $$
 It follows that the eigenvalues of $\rho_W$ consist of all possible products
 $$
 \lambda_{jk} = \frac{1}{j(j+1)\pi^2 k^2}
 $$
 with corresponding eigenfunctions 
$$
f_{\chi,j}(s)f_{\psi,k}(t).
$$
The expansion in Theorem 1 is then 
Parseval's identity with 
$$
Z_{jk}= \int_0^1 \int_0^1 \mathbb{W}(s,t) f_{\chi,j}(s)f_{\psi,k}(t) \, ds \, dt.
$$

\subsection{Numerical Work} \label{subsec:numericalNull}
The distribution of $\Wbar_\infty$ can be computed numerically in order to provide approximate, asymptotically valid, p-values.  
Our desired approximation to the p-value is 
$$
P(\Wbar_n>w_{\rm obs}) \approx P(\Wbar_\infty > w_{\rm obs})
$$
where $w_{\rm obs}$ is the value of $\Wbar_n$ observed in the data.   
Define 
$$
\lambda_{jk} = \frac{1}{\pi^2 j(j+1) k^2}.
$$
In practice, we truncate the infinite sum defining $\Wbar_\infty$,  retaining the terms with the largest  values of $\lambda_{jk}$, and replace the neglected terms by their expected value. 
So we write
\begin{align*}
\Wbar_\infty & = \Wbar_M+T_M \\
& = \sum_{jk\le M} \lambda_{jk} Z_{jk}^2 + \sum_{jk> M} \lambda_{jk} Z_{jk}^2 .
\end{align*}
We then approximate $T_M$ by its expected value:
$$
\mu_M \equiv \sum_{jk> M} \lambda_{jk} {\rm E}\left( Z_{jk}^2 \right) = \sum_{jk> M} \lambda_{jk}.
$$
Since the mean of  $\Wbar_\infty$ is 
$$
\sum_{j,k} \lambda_{jk} = \frac{1}{6}
$$
the mean of $T_M$ may be computed by
$$
\frac{1}{6} - \sum_{jk \le M} \lambda_{jk}.
$$
Our approximation becomes 
$$
P(\Wbar_n > w_{\rm obs}) \approx P( \Wbar_M +\mu_M > w_{\rm obs}).
$$
The latter quantity may now be computed by using numerical Fourier inversion following \cite{imhof1961}. The \texttt{R} package \texttt{CompQuadForm} (see \citealp{duchesne2010computing}) implements this computation in the function \texttt{imhof}; we use this software in our numerical work below.  

\begin{figure}[H]
\includegraphics[width=\textwidth]{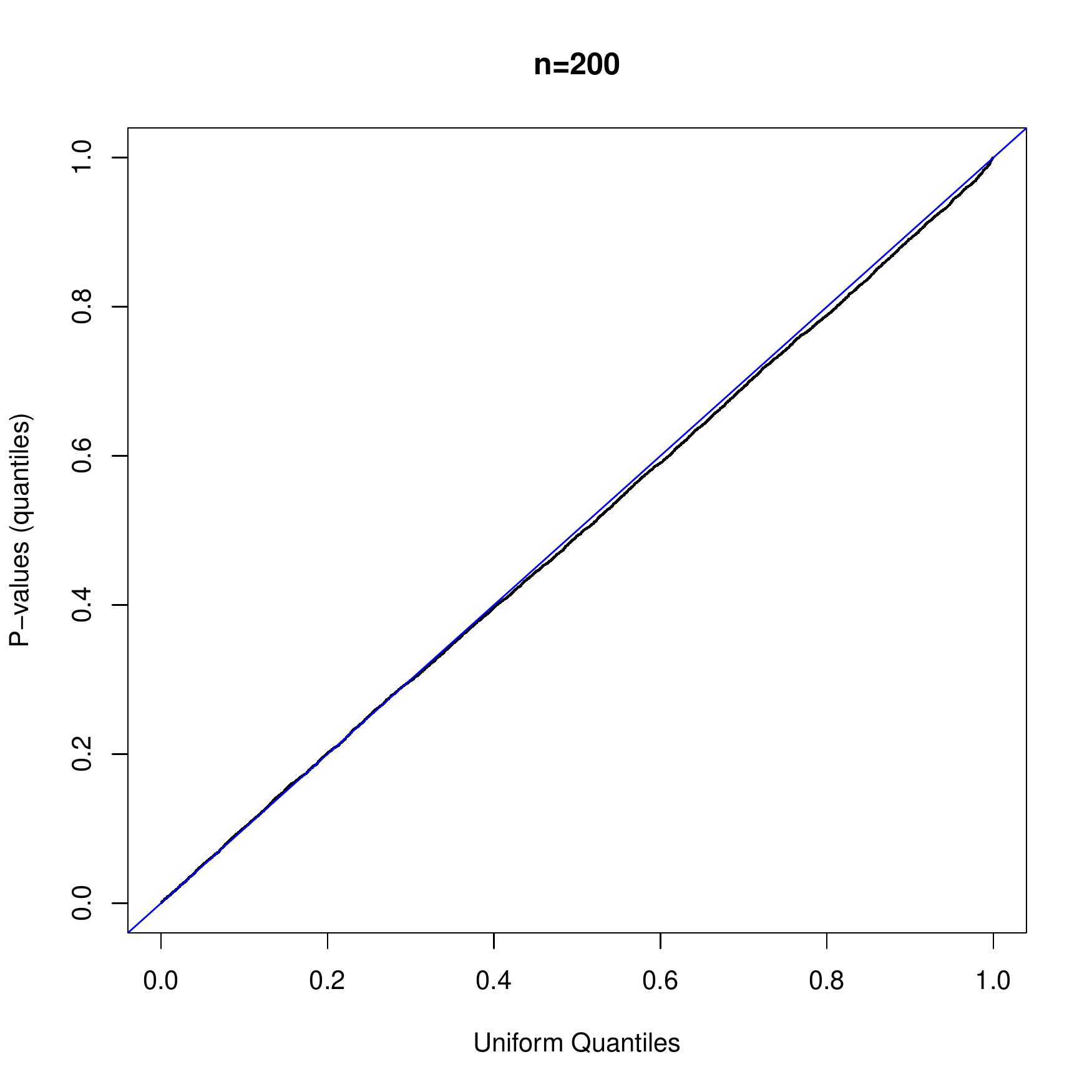}
\caption{Ordered p-values plotted against uniform quantiles for 10,000 iid Monte Carlo samples from a continuous distribution. The 
blue line is the uniform cumulative distribution function; exact p-values have a uniform distribution; the graph shows this approximation is good.}\label{fig:NullPvals200}
\end{figure}

\begin{figure}[H]
\includegraphics[width=\textwidth]{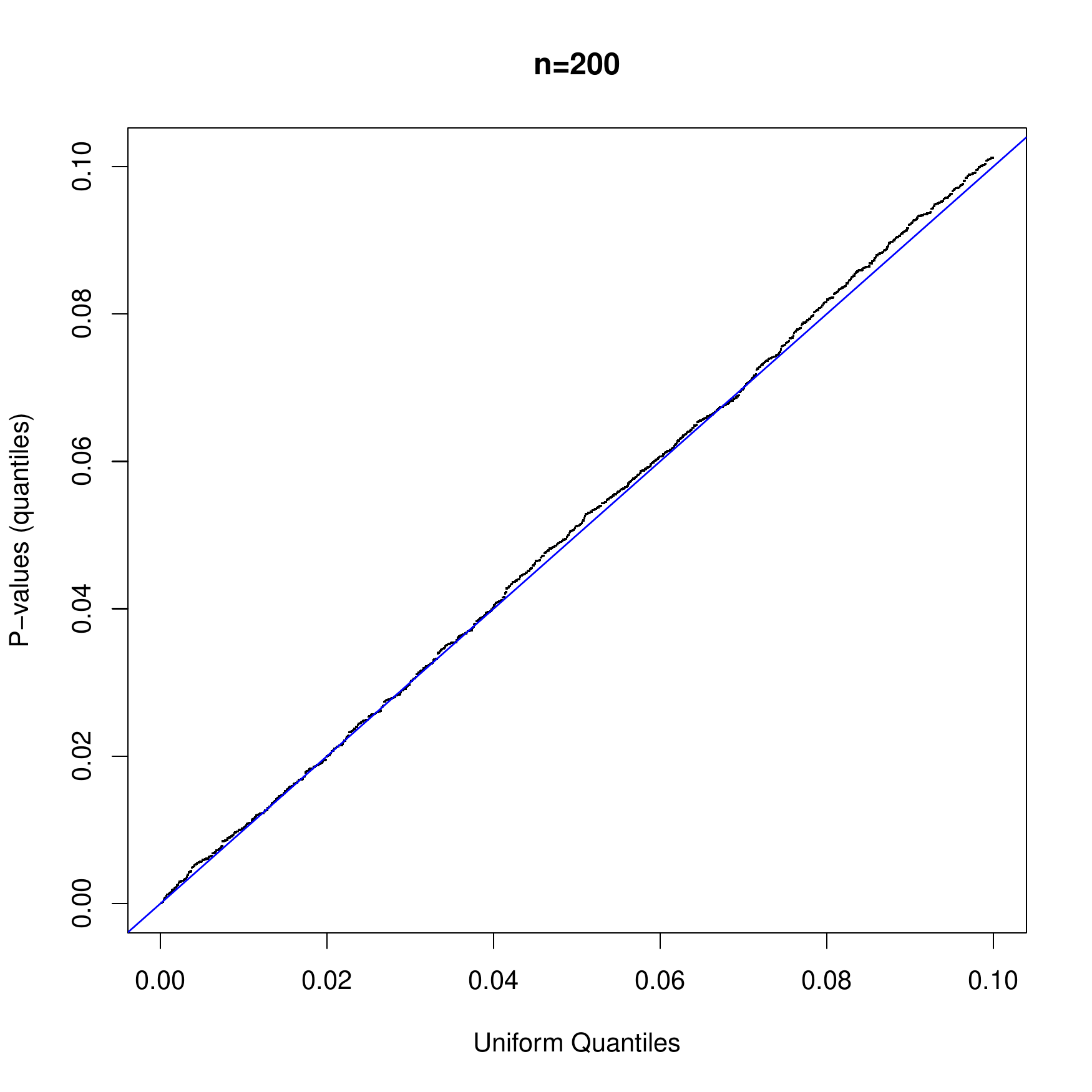}
\caption{Exploded view of Figure~\ref{fig:NullPvals200} showing the lower 10\% of the distribution of the ordered p-values plotted against uniform quantiles for 10,000 iid Monte Carlo samples from a continuous distribution.  The 
blue line is the uniform cumulative distribution function; exact p-values have a uniform distribution; the graph shows this approximation is very good in the important lower tail.}\label{fig:NullPvals200Zoom}
\end{figure}

We have evaluated the quality of our asymptotic approximation to the null distribution of $\Wbar$ in a small Monte Carlo study.  Since this distribution does not depend on $H$ when the null hypothesis holds we generated $N=10,000$ samples of size $n=200,500,1000$.  Figure~\ref{fig:NullPvals200}  shows a Q-Q plot for these 10,000 values  for $n=200$ to check the uniformity of their distribution. Specifically, we plot the order statistics against the uniform plotting points $1/(N+1),\ldots,N/(N+1)$.  Figure~\ref{fig:NullPvals200Zoom} is an enlargement of the smallest 10\% of these values since the quality of the approximation is most important for small p-values.  In both cases it is seen that the approximation is excellent.  For completeness, however, we note that the hypothesis of exact uniformity of these 10,000 p-values is rejected ($P\approx 0.01$) by the Anderson-Darling test.  Applied to the smallest 1,000 p-values, rescaled so that p-value number 1,001 from the bottom becomes 1, the Anderson-Darling p-value is actually 0.99.  We conclude the uniform approximation is very good at reasonable sample sizes, particularly in the important lower tail. For p-values over 0.5 we believe that the truncation we must do in order to compute the limit law is slightly off but argue that inaccuracy in the upper tail of p-values is not very consequential.

\section{Monte Carlo Power approximations} \label{sec:PowerMC}

We undertook a variety of Monte Carlo simulation studies to compare the power of $\overline{W}_n$ to $W_{\max}$. In Table \ref{table:power} we show the percentage of samples rejected in 10,000 trials by the two methods at the levels $\alpha=0.05$ and $\alpha=0.1$. We consider samples of size $n\in\{20,50,100\}$.  In one experiment recorded in the table we generated data from the Gamma distributions where the parameters change at $c=n/2$. In another experiment we change from the Gamma distribution to the Normal distribution at $c=n/2$; in this case neither the mean nor the variance changes. While our tests are designed to detect single change points we have included two trials in which  there are three segments which change between various Gamma distributions. One changes from shape 1, scale 2 to shape 2, scale 1 at the 40\% point and then to shape 0.5, scale 4 at the 60\% point. All three of these have the same mean.  The other changes from shape 1, scale 2 to shape 2, scale 3, and back to shape 1, scale 2; the changes happen after 30\% and then 70\% of the data.  Finally we present two experiments with samples from the normal distribution; in one the mean changes at $c=n/2$ and in the other the standard deviation changes at the same point.  In all these trials the parameter values in the distributions in a given segment do not change as the sample size changes; this may be compared with the further Monte Carlo results in Section~\ref{sec:ContiguousPower}. 

\begin{table}[H]\label{table:power}
	\begin{center}
		\resizebox{\columnwidth}{!}{%
			\begin{tabular}{ c | c | c | c | c | c }
				\multicolumn{2}{c}{} & \multicolumn{2}{c}{$\alpha=0.1$} & \multicolumn{2}{c}{$\alpha=0.05$}\\
				\hline 
				Alternative & Sample size & $W_{\text{max}}$ & $\overline{W}_n$ & $W_{\text{max}}$ & $\overline{W}_n$ \\
				\hline
				\multirow{3}{*}{\makecell{$X_1,\ldots ,X_{0.5n}\sim \text{Gamma}(1,2)$, \\ $X_{0.5n+1},\ldots ,X_{n}\sim \text{Gamma}(2,2)$}} & $n=20$ & $47.9$ & $50.7$ & $35.0$ & $37.5$\\
				& $n=50$ & $82.3$  & $85.7$ & $73.9$ & $77.4$ \\
				& $n=100$ & $98.3$ & $98.9$ & $96.3$ & $96.9$ \\
				\hline
				\multirow{3}{*}{\makecell{$X_1,\ldots ,X_{0.5n}\sim \text{Gamma}(1,2)$, \\ $X_{0.5n+1},\ldots ,X_{n}\sim \mathcal{N}(2,2)$}} & $n=20$ & $12.9$ & $13.7$ & $6.9$ & $7.2$ \\
				& $n=50$ & $16.1$  & $19.2$ & $9.0$ & $11.2$ \\
				& $n=100$ & $22.1$  & $31.2$ & $13.7$ & $19.0$ \\
				\hline
				\multirow{3}{*}{\makecell{$X_1,\ldots ,X_{0.4n}\sim \text{Gamma}(1,2)$, \\ $X_{0.4n+1},\ldots ,X_{0.6n}\sim \text{Gamma}(2,1)$ \\ $X_{0.6n+1},\ldots ,X_{n}\sim \text{Gamma}(0.5,4)$}} & $n=20$ & $17.5$ & $16.5$ & $10.0$ & $9.2$ \\
				& $n=50$ & $24.6$  & $25.5$ & $15.5$ & $15.9$ \\
				& $n=100$ & $38.3$  & $42.8$ & $27.3$ & $28.5$ \\
				
				\hline
				\multirow{3}{*}{\makecell{$X_1,\ldots ,X_{0.3n}\sim \text{Gamma}(1,2)$, \\ $X_{0.3n+1},\ldots ,X_{0.7n}\sim \text{Gamma}(2,3)$ \\ $X_{0.7n+1},\ldots ,X_{n}\sim \text{Gamma}(1,2)$}} & $n=20$ & $29.0$ & $20.6$ & $15.8$ & $7.9$  \\
				& $n=50$ & $72.3$  & $71.6$ & $54.4$ & $48.1$ \\
				& $n=100$ & $98.3$  & $98.6$ & $94.1$ & $94.6$ \\ 
				\hline
				\multirow{3}{*}{\makecell{$X_1,\ldots ,X_{0.5n}\sim \mathcal{N}(0,1)$, \\ $X_{0.5n+1},\ldots ,X_{n}\sim \mathcal{N}(0,3)$}} & $n=20$ & $18.2$ & $22.0$ & $10.8$ & $11.3$\\
				& $n=50$ & $29.6$  & $56.0$ & $17.0$ & $33.0$ \\
				& $n=100$ & $66.3$ & $93.4$ & $45.0$ & $81.2$ \\
				\hline
				\multirow{3}{*}{\makecell{$X_1,\ldots ,X_{0.5n}\sim \text{Exp}(1)$, \\ $X_{0.5n+1},\ldots ,X_{n}\sim \text{Exp}(1.5)$}} & $n=20$ & $15.8$ & $16.4$ & $9.1$ & $9.3$\\
				& $n=50$ & $23.4$  & $26.9$ & $14.9$ & $17.5$ \\
				& $n=100$ & $35.8$ & $42.7$ & $25.0$ & $31.0$ 
			\end{tabular}
		}
	\end{center}
	\caption{Powers (percentage) from various alternative distributions and significance levels $0.1$ and $0.05$. Critical points were calculated with $100,000$ and Powers by $10,000$ Monte Carlo simulations. The notation Gamma($\alpha,\beta$) indicates sampling from a Gamma distribution with shape $\alpha$ and scale $\beta$. The parameters in the normal distribution are mean and variance as usual. The parameter in the Exponential distribution is the mean.}
\end{table}

It will be seen that, except for very small samples, when there is a single change point the test using $\overline{W}_n$ has better power than $W_{\rm max}$. Since it is also far faster to compute p-values for $\overline{W}_n$ using the highly accurate asymptotic law we recommend $\overline{W}$ over $W_{\rm max}$.  At the same time we observe that the procedure is specifically designed to choose between 1 change point and no change points and not to estimate and find multiple change points. In particular, for one of the alternatives in Table~\ref{table:power} with 2 change points the statistic $W_{\rm max}$ is usually more sensitive than $\overline{W}_n$.  

The results presented here show how the powers grow with sample size when the two distributions are fixed.  Other experiments, not reported here, show that both statistics have better power when the change is near the center of the sequence. More Monte Carlo power calculations are presented in Section~\ref{sec:Wmax} below with a focus on contiguous alternatives.

\section{Power approximations: contiguous alternatives}\label{sec:ContiguousPower}

We now compute approximate distribution theory for $\bar{W}_n$ when the null hypothesis is false and the extent of the change at the change point is big enough to be detectable but not obvious; that is, we study situations where the best possible power in large samples stays away from 1.  To do so we consider a sequence of alternatives indexed by $n$ and assume that these alternatives are contiguous to a sequence for which the null hypothesis of no change holds.  To be specific our null hypothesis sequence will have $X_i$ iid for $1 \le i \le n$ with density $h$ and cdf $H$.  For the alternative we suppose that there is a value $c_0$ such that for $1 \le i \le c_0$, the $X_i$ are iid with density $f$ and that for $c_0+1 \le i \le n$ the $X_i$ are iid with density $g$. All of $f$, $g$, $h$, and the true change point $c_0$ may depend on $n$ but the dependence will be hidden in our notation.  Under the null hypothesis the joint density of $X_1,\ldots, X_n$ is 
$$
{\bf f}_{0 n}(x_1,\ldots,x_n) = \prod_{i=1}^n h(x_i).
$$
Under the alternative the joint density becomes
$$
{\bf f}_{1 n}(x_1,\ldots,x_n) = \prod_{i=1}^{c_0} f(x_i) \prod_{c_0+1}^n g(x_i).
$$
The log-likelihood ratio of these two is 
\begin{align*}
\Lambda_n &= \ln \left\{{\bf f}_{1,n}(X_1,\ldots,X_n)/{\bf f}_{0 n}(X_1,\ldots,X_n)\right\} 
\\
&= \sum_{i=1}^{c_0} \ln\left\{f(X_i)/h(X_i)\right\}+\sum_{i=c_0+1}^n \ln\left\{g(X_i)/h(X_i)\right\}.
\end{align*}
The sequence of alternatives ${\bf f}_{1 n}$ is contiguous to the null sequence ${\bf f}_{0 n}$ if, computing under the null hypothesis, we have
\begin{equation}\label{eq:contiguity}
\Lambda_n \rightsquigarrow N(-\tau^2/2, \tau^2)
\end{equation}
for some $0 \le \tau < \infty$.  If we define $U_i=H(X_i)$ then under the null hypothesis the $U_i$ are iid Uniform[0,1]. Under the alternative $U_1,\dots,U_{c_0}$ are iid with density
$\tilde{f}(u) = f(H^{-1}(u))/h\left(H^{-1}(u)\right)$ while $U_{c_0+1},\ldots,U_n$ are iid with density $\tilde{g}(u) = g(H^{-1}(u))/h\left(H^{-1}(u)\right)$. The likelihood ratio becomes
$$
\tilde\Lambda_n = \sum_{i=1}^{c_0} \ln\left\{\tilde{f}(U_i)\right\}+\sum_{i=c_0+1}^n \ln\left\{\tilde{g}(U_i)\right\}.
$$
Since our test statistics are invariant to a monotone transformation applied to each individual data point we will take $H$ to be Uniform[0,1] and then drop the tildes from our notation.
The quantity
$$
 S_n =  \sum_{i=1}^{c_0} \phi_f(X_i)/\sqrt{n} +\sum_{i=c_0+1}^n \phi_g(X_i)/\sqrt{n} 
$$
is needed in our theorem.

\begin{theorem} \label{theorem:contiguous} Assume
\begin{description}

\item[{\bf A1}] There are two functions $\phi_f$ and $\phi_g$ in $L_2[0,1]$ such that
$$
\lim_{n\to\infty} \sqrt{n}(f-1) = \phi_f
$$
and
$$
\lim_{n\to\infty} \sqrt{n}(g-1) = \phi_g.
$$
\item[{\bf A2}] There is a $u\in(0,1)$ such that
$$
\lim_{n\to\infty}\frac{c_0}{n} = u.
$$
\end{description}
Then as $n \to \infty$ we have, under the sequence of alternative hypotheses specified by $f$, $g$, and $c$,
\begin{enumerate}

\item 
The log-likelihood ratio satisfies
$$
\Lambda_n = S_n+o_P(1)  \rightsquigarrow N(-\tau^2/2,\tau^2)
$$
where 
$$
\tau^2 = u \int_0^1 \phi_f^2(t)\, dt +(1-u) \int_0^1 \phi_g^2(t)\, dt .
$$

\item The process $\mathbb{W}_n$ converges weakly to a Gaussian process with covariance $\rho$ and mean
$$
\mu(s,t)= \mu_\chi(s)\mu_\psi(t)
$$
where
$$
\mu_\chi(s) =  \sqrt{s(1-s)}\left\{\frac{1-u}{1-s} 1(s\le u) + \frac{u}{s} 1(s>u)\right\} 
$$
and
$$
\mu_\psi(t) =
\left[{\rm E}\left\{\phi_f(U)1(U \le t)\right\} - {\rm E}\left\{\phi_g(U)1(U \le t)\right\} \right].
$$

\item and
$$
\Wbar_n  \rightsquigarrow  \Wbar_\infty \equiv \sum_{j=1}^\infty \sum_{k=1}^\infty \frac{\left(Z_{jk}+\eta_j\tau_k\right)^2}{j(j+1) \pi^2 k^2}
$$
 where the  $Z_{jk}$ are iid standard normal,
 $$
 \eta_j =\int_0^1 \mu_\chi(s) f_{j,\chi}(s) \, ds,
 $$
 and
 $$
 \tau_k =  \int_0^1 \mu_\psi(t) f_{j,\psi}(t) \, dt.
 $$
 \end{enumerate}
 
 \end{theorem}

As with the null distribution, this limiting alternative distribution for $\Wbar$ can be computed using the \texttt{R} package \texttt{CompQuadForm}.
As an example we take $f$ to be standard normal and $g$ to be normal with mean $\mu$ and standard deviation $\sigma$. The  two parameters are assumed to depend on $n$ in such a way that 
$$
\sqrt{n} \mu \to \gamma_1 \quad \text{ and } \quad \sqrt{n}(\sigma-1) \to \gamma_2.
$$
It is convenient to take $h=f$.  Under the null  the data $X_1,\ldots,X_n$ are iid standard normal.  The functions $\tilde{f}$ and $\tilde{g}$ are then given by
$\tilde{f}\equiv 0$ and 
$$
\tilde{g}(u) = \frac{\phi\left\{\frac{\Phi^{-1}(u) - \mu}{\sigma}\right)}{\phi\left\{\Phi^{-1}(u)\right\}}.
$$
Under these conditions  we may check that condition \textbf{A1} holds with 
$\phi_f=0$ 
and 
$$
\phi_g(u) = \gamma_2\left[\left\{\Phi^{-1}(u)\right\}^2-1\right] + \gamma_1 \Phi^{-1}(u).
$$

\section{Large sample behaviour of $W_{\rm max}$}  \label{sec:Wmax}

The statistic $W_{\rm max}$ is more challenging to analyze because the weak convergence result in Theorem~\ref{theorem:weakconvergence} asserts convergence in 
$\ell_\infty^{\rm loc}((0,1) \times[0,1])$. By $ \ell_\infty^{\rm loc}((0,1) \times[0,1])$ we mean the space of functions on $(0,1)\times[0,1]$ which are bounded on 
compact subsets of their domain.  We give this the topology of uniform convergence on compacts.  See \cite{van1996weak}.  Our proof of Theorem~\ref{theorem:wbar} shows that our statistic is a continuous function on a subset of $\ell_\infty^{\rm loc}((0,1) \times[0,1])$ to which sample paths of $\mathbb{W}_\infty$ are almost sure to belong. We are not able to establish the corresponding result for $W_{\rm max}$.  Traditionally this problem has been handled either by fixing a small $\epsilon>0$ and redefining $W_{\rm max}$ by maximizing only over $\{c: \epsilon \le c/n \le 1-\epsilon\}$ or by careful analysis of the behaviour of the process and the test statistic for $c/n$ close to 0 or to 1. For instance, \cite{jaeschke1979} considers a weighted Kolmogorov-Smirnov test for the uniform distribution and shows that the supremum of the weighted empirical process has, after suitable normalization, an extreme value distribution.

\begin{figure}
\begin{center}
\includegraphics[width=\textwidth]{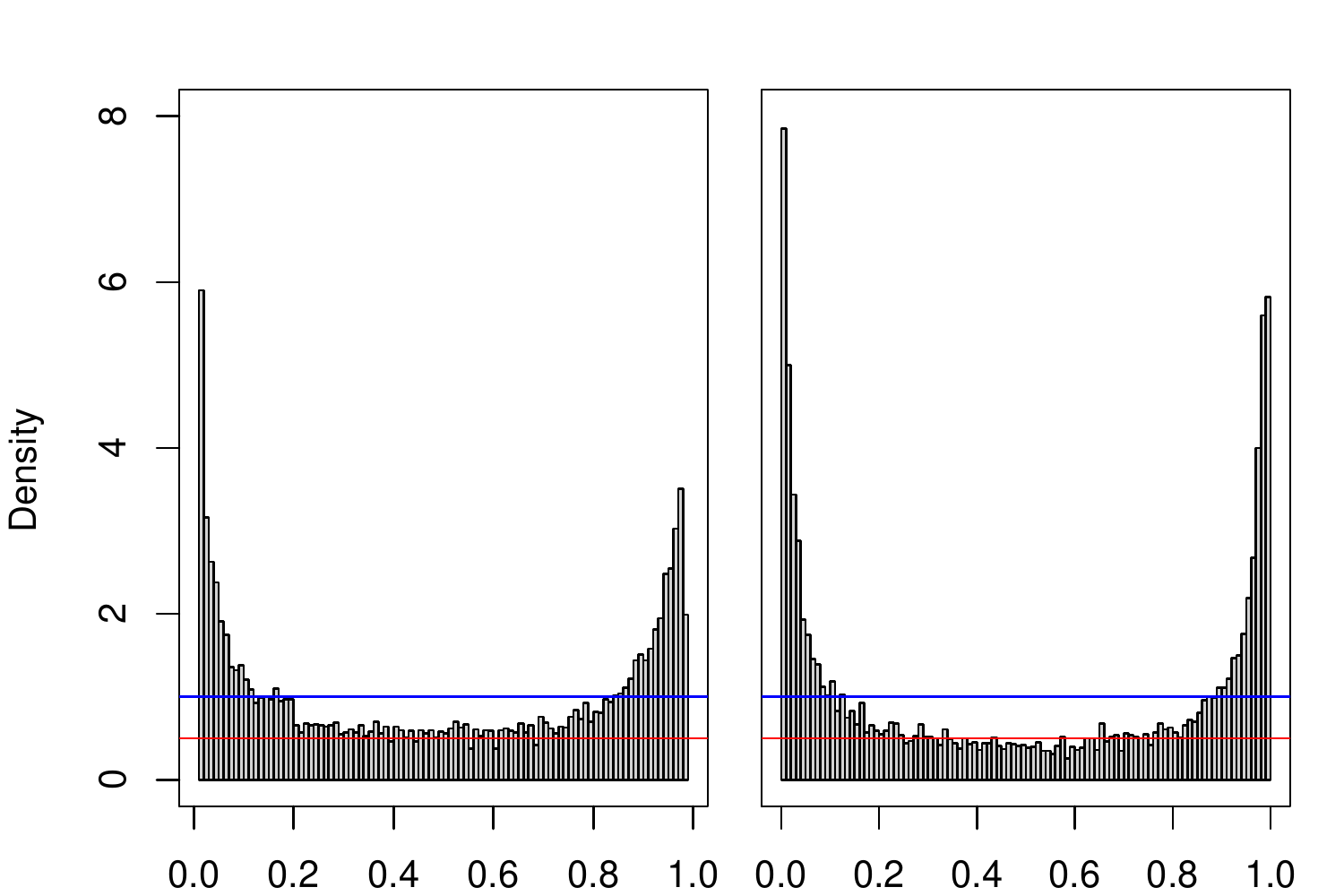}
\end{center}
\caption{Histograms of values of estimated change points for sample sizes $n=100$ on the left and $n=500$ on the right.  The null hypothesis is true and 10,000 samples were used for each histogram. The $x$-axis shows $\hat{c}/n$ and the $y$-axis is a probability density scale. The two figures have the same scales on each axis. Horizontal lines at height 1 (blue) and 0.5 (red) are provided to help see the extent to which the distribution on the right is more concentrated around 0 and 1 than the distribution on the left.}\label{fig:chat}
\end{figure}

We have not pursued either of these ideas but offer here some evidence that this statistic has some important defects.  First we look at a small simulation study.  We generated 10,000 samples of size 100 and 500 from the null hypothesis. In Figure~\ref{fig:chat} we plot histograms of the value $\hat{c}$ which maximizes $W_n(c)$ over $1 \le c \le n-1$.  Observe that as the sample size grows the histogram concentrates near 0 and 1 (though the convergence is slow). We can prove:

\begin{proposition}\label{prop:chatbehaviour}
Under the null hypothesis and under any sequence of contiguous alternatives 
$$
\min\left\{\frac{\hat{c}}{n}, \frac{n-\hat{c}}{n}\right\} \to 0
$$
in probability. Under the null hypothesis, the distribution of $\hat{c}/n$ converges to a Bernoulli$(0.5)$ law.
\end{proposition}

This means that, even for data from detectable (but not obvious) alternatives, our test statistic $W_{\rm max}$ usually compares the distribution of a tiny fraction of the data to that of the vast majority of the data even when the true change point is in the middle of the sequence.  We also conjecture:
\begin{conjecture} \label{conj:zeroARE}
For any sequence of contiguous alternatives the difference between the power and the level of a test based on $W_{\rm max}$ goes to 0 as $n\to\infty$.
\end{conjecture}

\begin{table}[H]\centering
\begin{tabular}{cccccccc}
\multicolumn{7}{c}{Gamma, shape=$1+b/\sqrt{n}$, break at $n/2$}\\\\
 & & & $n=10$ & $n=50$ & $n=100$ & $n=200$ & $n=500$\\
&$\bar{W}$ & MC & 11.70 & 13.96 & 14.83 & 14.71 & 15.91 \\
$b=2$ &$\bar{W}$ & Asym & 11.79 & 13.59 & 14.61 & 14.67 & 15.70 \\
&$W_{\rm  max}$ & MC & 12.13 & 12.00 & 12.36 & 11.41 & 11.80 \\
\hline

&$\bar{W}$ & MC & 18.50 & 25.18 & 26.48 & 27.74 & 29.52 \\
$b=3$ &$\bar{W}$ & Asym & 18.72 & 24.73 & 26.12 & 27.66 & 29.25 \\
&$W_{\rm  max}$ & MC & 18.62 & 22.05 & 21.84 & 21.34 & 21.88 \\
\hline
&$\bar{W}$ & MC & 34.95 & 52.67 & 57.39 & 61.28 & 65.62 \\
$b=5$ &$\bar{W}$ & Asym & 35.26 & 51.97 & 57.06 & 61.18 & 65.35 \\
&$W_{\rm  max}$ & MC & 35.48 & 47.60 & 50.07 & 52.76 & 54.46 \\
\hline
 \\
&\multicolumn{7}{c}{Gamma, shape=$1+b/\sqrt{n}$, break at $3n/10$}\\\\
&$\bar{W}$ & MC & 9.24 & 11.29 & 11.73 & 11.83 & 13.21 \\
$b=2$ & $\bar{W}$ & Asym &  9.42 & 10.86 & 11.47 & 11.79 & 13.10 \\
&$W_{\rm max}$ & MC & 10.00 & 10.60 & 10.48 & 9.86 & 10.37 \\
\hline
&$\bar{W}$ & MC & 13.41 & 20.04 & 20.26 & 21.80 & 23.15 \\
$b=3$ & $\bar{W}$ & Asym & 13.54 & 19.56 & 19.92 & 21.66 & 22.98 \\
&$W_{\rm  max}$ & MC & 14.81 & 18.07 & 17.38 & 18.00 & 17.97 \\
\hline
&$\bar{W}$ & MC & 22.42 & 41.53 & 45.54 & 48.59 & 53.34 \\
$b=5$ & $\bar{W}$ & Asym & 22.75 & 40.87 & 45.11 & 48.52 & 53.07 \\
&$W_{\rm  max}$ & MC & 26.43 & 39.44 & 41.36 & 43.09 & 45.72 \\
\hline \\
\end{tabular}
\caption{Powers (percentage) for change from Gamma(shape$=1+b/\sqrt{n}$, scale=1) to Gamma(1,1) at the indicated breakpoint, $n/2$ in the top and $3n/10$ in the bottom. Powers are based on 10,000 samples and either use Monte Carlo critical points (based on 100,000 samples) or asymptotic critical points as indicated by `MC' or `Asym'. All tests are at the level $\alpha=0.05$. }\label{table:gammacontiguouspower}
\end{table}

\begin{table}[H]\centering
\begin{tabular}{cccccccc}
\multicolumn{7}{c}{Normal, $\sigma =1+b/\sqrt{n}$, break at $n/2$}\\\\
& & & $n=10$ & $n=50$ & $n=100$ & $n=200$ & $n=500$\\
&$\bar{W}$ & MC & 5.61 & 5.97 & 5.65 & 5.66 & 5.91 \\
$b=2$ & $\bar{W}$ & Asym & 5.69 & 5.77 & 5.40 & 5.61 & 5.83 \\
&$W_{\rm  max}$ & MC & 6.70 & 5.72 & 5.19 & 4.80 & 5.25 \\
\hline \\[-11pt]
&$\bar{W}$ & MC & 6.11 & 7.04 & 6.87 & 6.75 & 7.40 \\
$b=3$ & $\bar{W}$ & Asym & 6.20 & 6.66 & 6.66 & 6.73 & 7.23 \\
&$W_{\rm max}$ & MC & 7.67 & 6.49 & 5.71 & 5.36 & 5.55 \\
\hline  \\[-11pt]
&$\bar{W}$ & MC & 6.76 & 9.55 & 11.10 & 11.32 & 13.56 \\
$b=5$ & $\bar{W}$ & Asym & 6.79 & 9.24 & 10.79 & 11.25 & 13.33 \\
&$W_{\rm max}$ & MC & 8.99 & 7.91 & 6.99 & 6.84 & 6.88 \\
\hline \\
&\multicolumn{7}{c}{Normal, $\sigma =1+b/\sqrt{n}$, break at $0.3n/10$}\\
\\
&$\bar{W}$ & MC & 6.26 & 6.49 & 5.80 & 5.63 & 5.76 \\
$b=2$ & $\bar{W}$ & Asym & 6.37 & 6.17 & 5.63 & 5.63 & 5.68 \\
&$W_{\rm max}$ & MC & 7.12 & 6.08 & 5.72 & 5.22 & 5.42 \\
\hline  \\[-11pt]
&$\bar{W}$ & MC & 6.91 & 7.37 & 6.74 & 6.41 & 6.95 \\
$b=3$ & $\bar{W}$ & Asym & 7.09 & 7.10 & 6.51 & 6.39 & 6.80 \\
&$W_{\rm max}$ & MC & 8.18 & 7.08 & 6.29 & 5.94 & 5.95 \\
\hline  \\[-11pt]
&$\bar{W}$ & MC & 7.89 & 9.40 & 9.92 & 9.91 & 11.13 \\
$b=5$ & $\bar{W}$ & Asym & 8.09 & 8.99 & 9.65 & 9.79 & 10.98 \\
&$W_{\rm max}$ & MC & 9.80 & 8.96 & 8.04 & 7.67 & 7.19 \\
\hline \\
\end{tabular}

\caption{Powers (percentage) for change from Normal(0,$\sigma=1+b/\sqrt{n}$) to Normal(0,1) at the indicated breakpoint, namely, $n/2$ in the top and $3n/10$ in the bottom. Powers are based on 10,000 samples and either use Monte Carlo critical points (based on 100,000 samples) or asymptotic critical points as indicated by `MC' or `Asym'. All tests are at the level $\alpha=0.05$. }\label{table:normalcontiguouspower}
\end{table}

Here is some Monte Carlo evidence from a simulation study. In Tables~\ref{table:gammacontiguouspower} and~\ref{table:normalcontiguouspower} we study four alternatives at sample sizes $n=10,50,100,200,500$.  For each sample size we draw 10,000 samples of size $n$. The first $c$ observations in each sample have some parameter of the form $a+b/\sqrt{n}$ and the remaining $n-c$ have parameter $a$. We used the Gamma distribution and the normal distribution and tried $c=0.5n$ and $c=0.3n$ for each distribution.  In the Gamma case we tried changing the shape parameter with $a=1$ while holding the scale parameter at 1.  The tables show the expected convergence (although we have not computed the power predicted by our theory in Section \ref{sec:ContiguousPower}. 

For the statistic $W_{\rm max}$ the tables  show, in the normal case, the power   declining towards the level (which is 5\% here).  For the Gamma cases studied here the power is rising but slowly for distant alternatives (large values of $b$) and declining very slowly for less distant alternatives (smaller values of $b$). Our experience in general is that for more distant alternatives it requires larger sample sizes before the power of $W_{\rm max}$ begins to drop.

Our conjecture is motivated by an analogy with \cite{lockhart1991overweight} in which it is shown that goodness-of-fit test statistics which depend only on $o(n)$ tail order statistics have the property asserted in the second conjecture. In the Appendix we prove the proposition and provide partial details showing how we would hope to prove our conjecture, if we could.

\section{Discussion}\label{sec:Discussion}
 
 It is a general principle that procedures with optimal frequency properties are found by searching among Bayes procedures.  It is also generally the case that optimal Bayes procedures involve averaging rather than maximizing.  These heuristics motivate considering testing for change points by using test statistics which are averages over possible change points rather than maxima. In this paper we have used this heuristic to motivate an average two sample goodness of fit statistic when we are concerned about general changes in distribution, rather than simple changes in mean, in a sequence of independent data points.  We have shown the resulting test statistic has computable large sample theory which can be used to provide very accurate p-values. Moreover we have shown that averaging over possible change points is generally more sensitive to alternatives than maximizing over possible change points.
 
The basic idea can be used in other contexts. Consider, for instance, testing for a change in mean.  We describe first the unrealistic situation in which the standard deviation is known and then how to handle estimation of that SD. Suppose $X_1,\ldots,X_n$ are independent and we wish to test the null hypothesis that they are iid with unknown mean $\mu$ and known standard deviation $\sigma$ (which we take to be 1 for notational convenience) against the alternative that the mean changes after the data point number $c$.  The usual $Z$ statistic is 
$$
T_c = \left(\frac{X_1+\cdots + X_c}{c} - \frac{X_{c+1}+\cdots +X_n}{n-c}\right)/\sqrt{\frac{1}{c}+\frac{1}{n-c}}.
$$
Our proposal would be to use the two sided test
$$
\overline{T^2} = \frac{1}{n-1} \sum_{c=1}^n T_c^2.
$$
This statistic has mean 1 under the null hypothesis of no change in mean.  Arguments similar to those in Section~\ref{sec:NullLimit} show that this statistic has the same limiting distribution, under the null,  as the well known Anderson-Darling goodness-of-fit statistic.  

In the more reasonable case where the (assumed common) standard deviation is unknown will use the statistic
$$
T_s^2 = \overline{T^2} /s^2
$$
where $s^2$ is some estimate of $\sigma^2$ which is consistent under the null hypothesis. The sample standard deviation is one possibility though this can be badly biased under the alternative.  An estimate which is rather less precise but still likely to be quite accurate under the alternative hypothesis is 
$$
s_1^2 = \frac{\sum_{i=1}^{n-1} (X_{i+1} - X_i)^2}{2(n-1)}.
$$
Notice that under the alternative hypothesis all but one term in this average is an unbiased estimate of $\sigma^2$; the bias in the estimator is $\Delta_\mu^2/(2n)$ where $\Delta_\mu$ denotes the change in the mean at the true change point.  Under the null our estimate is unbiased. The statistic $T_s^2$ also has the same limiting distribution as the well known Anderson-Darling goodness-of-fit statistic when the null holds.

Other nonparametric goodness of fit tests can be used instead of the Cram\'er-von Mises test. For example a Bayesian test \cite{labadi2014twosample}, likelihood tests \cite{LimitTheorems} or other two-sample tests \cite{Robustpower}.  Sample size, the kind of alternative distribution from which we expect the data to come  and the expected index of the change point should likely be used to choose the best test. Finding the asymptotic distribution for less well-known tests can be difficult. Bootstrapping can be used instead. This deserves further research.

\section*{Appendix}

\noindent\textbf{Proof of Theorems~\ref{theorem:wbar}~and~\ref{theorem:weakconvergence}}.

The weak limit $\mathbb{Z}$ given below is discussed in~\cite{TimeSeriesChangePoint} but we provide details for completeness.

We prove Theorem~\ref{theorem:weakconvergence} first.  Define the partial sum empirical process \cite[p. 225]{van1996weak}, for $(s,t) \in [0,1]^2$, by
$$
\mathbb{Z}_n(s,t) = \frac{1}{\sqrt{n}} \sum_{1 \le i \le ns} \left\{1(U_i \le t) - t\right\}.
$$
Our statistic can be described in terms of this process. 
Notice that 
$$
F_c(t) = \frac{\sqrt{n}}{c} \mathbb{Z}_n(c/n,t)  + t
$$
and that
$$
G_d(t) = \frac{\sqrt{n}}{d} \left\{\mathbb{Z}_n(1,t)-\mathbb{Z}_n(c/n,t)\right\} + t.
$$
Thus 
$$
F_c(t) -G_d(t) = \sqrt{n}\left\{\frac{ \mathbb{Z}_n(c/n,t)}{c} -\frac{\mathbb{Z}_n(1,t)-\mathbb{Z}_n(c/n,t)}{d}\right\}.
$$
Now define the process $\mathbb{W}_n(s,t)$ for $0<s<1$ and $0 \le t \le 1$ by
$$
\mathbb{W}_n(s,t) = \sqrt{s(1-s)} \left\{ \frac{\mathbb{Z}_n(s,t)}{s} - \frac{\mathbb{Z}_n(1,t)-\mathbb{Z}_n(s,t)}{1-s}\right\}.
$$
For given $c$ our  two sample test statistic is given by 
$$
W_n(c) = \int_0^1 \left\{\mathbb{W}_n(c/n,t)\right\}^2\, dH_n(t).
$$
Let $\nu_n$ be the probability measure on $(0,1)$ putting mass $1/(n-1)$ on each point of the form $c/n$ for $1 \le c \le n-1$.  Our statistic is
$$
\Wbar_n = \int_0^1\int_0^1 \left\{\mathbb{W}_n(s,t)\right\}^2\, dH_n(t)\, d\nu_n(s).
$$

We now break the proof of our two results into steps consisting of  a statement followed by a detailed proof. In each case the assertions are intended to hold under the null hypothesis and the assumption that the common distribution $H$ is continuous.

\medskip

\noindent\textit{Step 1}:  The process $\mathbb{Z}_n$ converges weakly in $\ell_\infty([0,1]^2)$ to a 
 tight, centred, Gaussian process $\mathbb{Z}$ with covariance 
$$
\rho(s,t;s',t') = (s\wedge s')(t \wedge t' - tt').
$$
See \cite{van1996weak}.  

\medskip

\noindent\textit{Step 2}:  Hence the process $\mathbb{W}_n$ converges weakly in $\ell_\infty^{\rm loc}((0,1) \times[0,1])$ to the tight centred Gaussian process
$$
\mathbb{W}(s,t) = \sqrt{s(1-s)} \left\{ \frac{\mathbb{Z}(s,t)}{s} - \frac{\mathbb{Z}(1,t)-\mathbb{Z}(s,t)}{1-s}\right\}.
$$
This process has continuous sample paths (on $(0,1) \times[0,1]$) and the covariance given in the statement of the theorem.

\medskip

\noindent\textit{Step 3}: 
For any sequence $c_n$ with $\epsilon_n\equiv c_n/n \to 0$ we have 
$$
\left\{\int_0^{\epsilon_n} + \int _{1-\epsilon_n}^1\right\} \left\{\mathbb{W}_n(c/n,t)\right\}^2\, dH_n(t) d\nu_n(s)= \frac{\sum_{i=1}^{c_n} W_n(i) + \sum_{i=n+1-c_n}^n W_n(i)}{n-1} \to 0
$$
in probability.  Under the null hypothesis the mean of $W_n(c)$ is $1/6+1/(6n)$; see \cite{anderson1962}.  The expected value of the indicated quantity is thus
$$
\frac{2c_n}{n-1} \left(\frac{1}{6}+\frac{1}{6n}\right) \to 0.
$$

\medskip

\noindent\textit{Step 4}:  The integral
$$
W_\infty = \int_0^1 \int_0^1 \mathbb{W}^2(s,t)dt\, ds
$$
is almost surely finite. Since all the variates involved are non-negative we may compute
$$
{\rm E}(W_\infty) ={\rm E}\left( \int_0^1 \int_0^1 \mathbb{W}^2(s,t)dt\, ds\right) = \int_0^1 \chi(s,s)\,ds \int_0^1 \psi(t,t)\, dt = 1/6<\infty.
$$

\medskip

\noindent\textit{Step 5}: 
For any sequence $\epsilon_n$ tending to 0 as $n\to\infty$ we have, by taking expectations,
$$
\left\{\int_0^{\epsilon_n} + \int_{1-\epsilon_n}^1\right\} \int_0^1 \mathbb{W}^2(s,t)dt\, ds \to 0
$$
in probability.

\medskip

\noindent\textit{Step 6}:
The tensor product kernel 
$$
\rho=\chi \otimes \psi (s,t;s',t') = \chi(s,s')\psi(t,t')
$$
is compact and has eigenvalue-eigenfunction pairs
$$
\lambda_{jk} =\frac{1}{j(j+1)}\frac{1}{\pi^2 k^2}, \quad f_{jk}(s,t) =f_{\chi,j}(s)f_{\psi,k}(t)
 $$
 indexed by $j,k$ each running from 1 to $\infty$.  It follows as usual that
 the family
 $$
 Z_{jk} = \frac{1}{\sqrt{\lambda_{jk}}} \int_0^1\int_0^1 \mathbb{W}(s,t) f_{jk}(s,t) dt\, ds
 $$
 defines a family of independent standard normal variables. Parseval's identity is then
 $$
 \int_0^1\int_0^1 \mathbb{W}^2(s,t)dt\, ds = \sum_{j=1}^\infty \sum_{k=1}^\infty \frac{Z_{jk}^2}{j(j+1) \pi^2 k^2}.
 $$
 
\medskip

 \noindent\textit{Step 7}:
 For each fixed $\epsilon>0$ we have 
 $$
 \int_\epsilon^{1-\epsilon} \int_0^1 \mathbb{W}_n^2(s,t) dH_n(t) \, d\nu_n(s) - \frac{1}{n-1} \sum_{n\epsilon < i < n(1-\epsilon)} W_n^2(i) \to 0
 $$
 in probability. This is an easy consequence of the fact that for $ i/n \le s < (i+1)/n$ we have $ \int_0^1\mathbb{W}_n^2(s,t) dF_n(t) = \mathbb{W}_n^2(i)$.
 
\medskip

 \noindent\textit{Step 8}:
 For each fixed $\epsilon>0$ we have 
 $$
 \int_\epsilon^{1-\epsilon} \int_0^1 \mathbb{W}_n^2(s,t) dH_n(t) \, d\nu_n(s) -  \int_\epsilon^{1-\epsilon} \int_0^1 \mathbb{W}_n^2(s,t) \,dt \, ds \to 0
 $$
Under the null hypothesis $H_n$ converges weakly to the uniform law on the unit interval.  Moreover $\nu_n$ converges weakly to Lebesgue measure on the unit interval.  The weak convergence result in Step 2 above uses a topology of uniform convergence on compacts such as the set $[\epsilon,1-\epsilon]\times[0,1]$ and this implies the desired result.
 
\medskip

 \noindent{\textit{Step 9}}: For each fixed $\epsilon>0$ we have 
 $$
 \int_\epsilon^{1-\epsilon} \int_0^1 \mathbb{W}_n^2(s,t) dt \, ds  \rightsquigarrow  \int_\epsilon^{1-\epsilon} \int_0^1 \mathbb{W}^2(s,t) \,dt \, ds .
 $$
This is a direct consequence of weak convergence using the continuous mapping theorem.

\medskip

 \noindent{\textit{Step 10}}: There is a metric $d$ on the set of probability measures on the real line for which the metric topology is the topology of weak convergence.  For each fixed $\epsilon>0$ we have 
$$
d\left({\mathcal L}\left( \int_\epsilon^{1-\epsilon} \int_0^1 \mathbb{W}_n^2(s,t) dt \, ds\right), {\mathcal L}\left(\int_\epsilon^{1-\epsilon} \int_0^1 \mathbb{W}^2(s,t) \,dt \, ds\right)\right) \to 0.
$$
There is then a sequence $\epsilon_n \to 0$ so slowly that this convergence continues to hold with $\epsilon$ replaced by $\epsilon_n$ and so that the convergences in Steps to 7 and 8 continue to hold.  Notice that by Step 5
$$
d\left({\mathcal L}\left(\int_{\epsilon_n}^{1-\epsilon_n} \int_0^1 \mathbb{W}^2(s,t) \,dt \, ds\right),{\mathcal L}\left(\int_0^{1} \int_0^1 \mathbb{W}^2(s,t) \,dt \, ds\right)\right)\to 0.
$$
for this sequence.
 
\medskip

  \noindent{\textit{Step 11}}:  For the sequence chosen in Step 10 we therefore have 
  $$
 \frac{1}{n-1} \sum_{n\epsilon_n < i < n(1-\epsilon_n)} W_n^2(i)  \rightsquigarrow  \int_0^{1} \int_0^1 \mathbb{W}^2(s,t) \,dt \, ds.
 $$
 In view of Step 1 we see 
 $$
 W_n  \rightsquigarrow  \int_0^{1} \int_0^1 \mathbb{W}^2(s,t) \,dt \, ds
 $$
 The law of the limit is, by Step 6, that of 
 $$
  \sum_{j=1}^\infty \sum_{k=1}^\infty \frac{Z_{jk}^2}{j(j+1) \pi^2 k^2}.
 $$
 This completes the proofs of Theorems~\ref{theorem:wbar} and~\ref{theorem:weakconvergence}.
 
 \bigskip
 
  \noindent{\bf Proof of Theorem~\ref{theorem:contiguous}}.
  
This is standard so we present only an outline.   
 Conditions \textbf{A1} and \textbf{A2} can be used to prove that
 $$
 \Lambda_n  - S_n  \to 0
 $$
 in probability under the null.  The Lindeberg Central limit theorem then establishes the first conclusion of the Theorem.  For more detailed arguments in a similar context see  \cite{guttorp1988asymptotic}. Thus, under the conditions of the theorem the sequence of alternatives is contiguous to a sequence for which the null holds.  
 
Contiguity implies that tightness under the null sequence extends to tightness under the alternative sequence.
This proves tightness, under the alternative, of the sequence  of processes $\mathbb{W}_n$.  Thus we need only compute the limiting finite dimensional distributions under the alternative sequence.  As usual we apply LeCam's Third Lemma (again similar arguments are in \cite{guttorp1988asymptotic})  to reduce the problem to studying the joint law, under the null hypothesis,  of $\Lambda_n$ and the vector $(\mathbb{W}_n(s_1,t_1),\ldots,\mathbb{W}_n(s_k,t_k)$ for an arbitrary sequence of time points $t_1,\ldots,t_k$ all in $[0,1]$.

The null distribution theory presented above (see Step 1 in the proof of Theorem~\ref{theorem:weakconvergence})
shows that, under the null hypothesis,
$$
\left(\mathbb{W}_n(s_1,t_1),\ldots,\mathbb{W}_n(s_k,t_k) \right)  \rightsquigarrow  MVN_k(0, \textbf{R}_W)
$$
where ${\bf R}_W$ is the $k\times k$ matrix with $i,j$th entry
$$
R_{Wij} = \rho_W(s_i,t_i; s_j,t_j).
$$
The Lindeberg Central Limit Theorem may now be used to show that the vector
 $$
 \left(S_n,\mathbb{W}_n(s_1,t_1),\ldots,\mathbb{W}_n(s_k,t_k)\right)
 $$
 converges in distribution to multivariate normal with mean vector $(-\tau^2/2,0,\ldots,0)$ and variance covariance matrix of the form
 $$
 \left[\begin{array}{cc}
  \tau^2 & \mathbf{c}^\top \\
  \mathbf{c} & \mathbf{R}_W
  \end{array}\right].
  $$
 Here the vector ${\bf c}$ is the limiting covariance which is found, after some algebra, to be
 $$
  c_i = \mu(s_i,t_i) = \mu_\chi(s_i)\mu_\psi(t_i).
  $$
  This completes the proof of the second assertion of the Theorem.
  
  The third step is standard; \cite{guttorp1988asymptotic} does similar problems.

\subsection*{Proof of  Proposition~\ref{prop:chatbehaviour}}

Fix $0 < \delta < 1/2$ and let $A_n$ denote the event $\{\delta \le \hat{c}/n \le 1-\delta\}$.  We will show that
$$
\lim_{n\to \infty} P(A_n) =0.
$$
This will prove Proposition~\ref{prop:chatbehaviour}.
To this end fix $0< \epsilon < \delta$.  Define
$$
M_n = \sup_{\delta \le s \le 1-\delta}
\int_0^1 \frac{\mathbb{B}_n^2(s,t)}{s(1-s)}\, dt
$$
and
$$
M_n'(\epsilon) = 
\sup_{\epsilon \le s \le \delta} \int_0^1\frac{\mathbb{B}_n^2(s,t)}{s(1-s)}\, dt.
$$
Then 
$$
A_n \subset \{M_n'(\epsilon) < M_n\}.
$$
Weak convergence of $\mathbb{B}_n$ to $\mathbb{B}$ guarantees that 
$$
\limsup_{n\to\infty} P(A_n)
\le \limsup_{n\to\infty} P\{M_n'(\epsilon) < M_n\}
\le P(M'(\epsilon)\le M)
$$
where 
$$
M =\sup_{\delta \le s \le 1-\delta}
\int_0^1\frac{\mathbb{B}^2(s,t)}{s(1-s)}\, dt
$$
and 
$$
M'(\epsilon) = 
\sup_{\epsilon < s \le \delta} \int_0^1\frac{\mathbb{B}^2(s,t)}{s(1-s)}\, dt.
$$ 
We claim that
\begin{equation}\label{eq:limsupbound}
\lim_{\epsilon\to 0}P(M'(\epsilon)\le M) = 0
\end{equation}
and this will prove  
$$
\limsup_{n\to\infty} P(A_n)
=0
$$
and Proposition~\ref{prop:chatbehaviour}.

Assertion~(\ref{eq:limsupbound}) would follow from a law of the iterated logarithm (as $s \to 0$) for the process
$$
s \mapsto \frac{\int_0^1 \mathbb{B}^2(s,t)\, dt}{s(1-s)}.
$$
While we expect such a result to hold we have not tried to prove anything along those lines. We will establish instead the lower bound 
$$
\limsup_{s\to 0}\int_0^1\frac{\pi^2\mathbb{B}^2(s,t)}{2\log\{\log(1/s)\}s(1-s)}\, dt \ge 1
$$
almost surely which is enough to imply~(\ref{eq:limsupbound}). We enumerate the steps needed:

\begin{enumerate}
    \item  Let 
    $$
    I_\mathbb{B}(s) = \int_0^1\frac{\mathbb{B}^2(s,t)}{s} \, dt,
    $$
    and
    $$
    I_\mathbb{Z}(s) = \int_0^1\frac{\mathbb{Z}^2(s,t)}{s} \, dt.
    $$
    Then 
\begin{align*}
   \int_0^1\frac{\mathbb{B}^2(s,t)}{s(1-s)}\, dt I &= \int_0^1\frac{\left\{\mathbb{Z}(s,t)-s\mathbb{Z}(1,t)\right\}^2}{s(1-s)}\, dt \\
   & \ge \frac{ I_\mathbb{Z}(s)}{1-s} + \frac{s}{1-s} \int_0^1\mathbb{Z}^2(1,t)\, dt -\frac{2s}{1-s} \sqrt{I_\mathbb{Z}(s)\int_0^1\mathbb{Z}^2(1,t) \, dt }.
   \end{align*}
   From this  we deduce that it is enough to show that 
\begin{equation} \label{eq:limsupZ}
  \limsup_{s\to 0}\int_0^1\frac{\pi^2\mathbb{Z}^2(s,t)}{2\log\{\log(1/s)\}s}\, dt \ge 1  
\end{equation}
almost surely.

\item For each fixed $s$ the process
$$
t \mapsto \frac{\mathbb{Z}(s,t)}{\sqrt{s}}
$$
is a Brownian Bridge. If we put
$$
W(s) = \int_0^1\frac{\mathbb{Z}^2(s,t)}{s}\, dt 
$$
then each $W(s)$ has the same distribution as the limit law of the usual Cram\'er-von Mises statistic which is the law of 
$$
\sum_{j=1}^\infty \lambda_j Z_j^2.
$$
In this representation the $Z_j$ are iid standard normal and the eigenvalues $\lambda_j$ are given, for $j=1,2,\ldots$, by
$$
\lambda_j = \frac{1}{\pi^2 j^2}.
$$

\item The process $\mathbb{Z}$ has independent increments in $s$ and  for each $0 < s' < s$ the process
$$
t\mapsto \frac{\mathbb{Z}(s,t)-\mathbb{Z}(s',t)}{\sqrt{s-s'}}
$$
has the same law as
$$
t \mapsto \frac{\mathbb{Z}(s,t)}{\sqrt{s}}
$$

\item Now fix $s_0=1$ and some $r<1$ to be chosen later. Define $s_n = s_0r^n$ for $n=1,2,\ldots$.
Put
$$
W_n = \int_0^1\frac{\mathbb{Z}^2(s_n,t)}{s_n}\, dt 
$$
and
$$
W_n^* = \int_0^1\frac{\left\{\mathbb{Z}(s_n,t)-\mathbb{Z}(s_{n+1},t)\right\}^2}{s_n-s_{n+1}}\, dt .
$$
All of these variables have the law of $W(s)$ described above.

\item  Fix $\epsilon>0$. Let $A_n$ be the event $W_n^* >2(1-\epsilon)\lambda_1\log(\log(1/s_n)) $ and $B_n$ be the event $W_{n+1} \le 2(1+\epsilon)\lambda_1\log(\log(1/s_n))$.  We will show that we can choose $r$ small enough so that
\begin{enumerate}
    \item The event that $A_n $ occurs infinitely often (i.o.) has probability 1.
    \item The event that $B_n$ occurs
    for all large $n$ has probability  1.
    \end{enumerate}
    
    \item So the event $A_n\cap B_n$ i.o. has probability 1.
    \item On the event $A_n\cap B_n$ we have
    $W_n \ge 2(1-\epsilon)\lambda_1\log(\log(1/s_n)$
    so that this event occurs infinitely often.
    
    \item This proves 
    $$
P\left\{W_n \ge 2(1-\epsilon)\lambda_1\log(\log(1/s_n))\, \text{i.o.}\right\}=1.
$$
which establishes~(\ref{eq:limsupZ}).
The definition of contiguity is that any sequence of events whose probability converges to 0 under the null has probability converging to 0 under the alternative. This finishes the proof of Proposition~\ref{prop:chatbehaviour}.

\end{enumerate}

\subsection*{Evidence for Conjecture 1}

For $\epsilon>0$ we define
$$
I_n(\epsilon) = \{c: 1 \le c \le n\epsilon \text{ or } 1 \le n-c \le n\epsilon \}.
$$
Proposition~\ref{prop:chatbehaviour} establishes that there is a sequence $\epsilon_n \searrow  0$ such 
$$
\lim_{n\to\infty} P(\hat{c}_n \in I_n(\epsilon_n)) = 1.
$$
Thus
\begin{equation} \label{eq:maxisatedges}
  P\left[W_{\rm max} = \max\{W_n(c): c\in I_n(\epsilon_n)\}\right] \to 1 .  
\end{equation}

We now outline the steps in our strategy for proving the conjecture before giving some evidence for each step.

\begin{itemize}

    \item[\textbf{Step 1}] There  are constants $a_n$ and $b_n$  and a random variable $V$ such that 
    $$
    a_n W_{\rm max} -b_n \rightsquigarrow V
    $$
    and $V$ has a continuous limit distribution. 
    
    \item[\textbf{Step 2}:] So
    $$
    a_n \max\{W_n(c): c\in I_n(\epsilon_n)\} -b_n \rightsquigarrow V.
    $$
    
    \item[\textbf{Step 3}:] There are random variables 
    $\tilde{W}_n(c)$ such that under the null hypothesis
    $$
    a_n \max\{|W_n(c)-\tilde{W}_n(c)|: c \in I_n(\epsilon_n)\}\to 0
    $$
    and such that for each $c\in I_n(\epsilon_n)$ the variable 
    $\tilde{W}_n(c)$ is measurable with respect to the $\sigma$field generated by $X_c, c\in I_n(\epsilon_n)$.
    To be specific we define, for $c<n/2$,
    $$
    \tilde{W}_n(c) = \int_0^1 c\{F_c(u)-u\}^2 \, du
    $$
    and, for $c>n/2$,
    $$
    \tilde{W}_n(c) = \int_0^1 d\{G_d(u)-u\}^2 \, du.
    $$
    (Recall the shorthand $d=n-c$.)
    
    \item[\textbf{Step 4}:] Define
    $$
    \tilde{\Lambda}_n = \sum_{n\epsilon_n < c \le c_0} \phi_f(X_c)/ \sqrt{n}  -
    \sum_{c_0 < c < n-n\epsilon_n} \phi_g(X_c)/\sqrt{n}.
    $$
    The log-likelihood ratio $\Lambda_n$ satisfies 
    $$
    \Lambda_n - \tilde{\Lambda}_n \to 0
    $$
    in probability, under the null hypothesis.
    
    \item[\textbf{Step 5}:] Since $\tilde{W}_{\rm max}$ is independent of  $\tilde{\Lambda}_n$ we may apply LeCam's third lemma to show that under the sequence of contiguous alternatives we have
    $$
    a_n W_{\rm max}-b_n \rightsquigarrow V
    $$
    
    \item[\textbf{Step 6}:] Since this limit law is the same as under the null we must power minus level tends to 0.

\end{itemize}

For some of these steps we can fill in partial evidence.

For Step 1 we would hope to follow the ideas in \cite{jaeschke1979} to show that the limit $V$ has an extreme value distribution. In that paper the maximizer of the usual empirical process, standardized by dividing by its standard deviation, is shown to have an extreme value limit with constants analogous to $a_n$ and $b_n$ involving $\sqrt{\log\log n}$ and $\log \log \log n$.

Step 2 is a consequence of Step 1 and~(\ref{eq:maxisatedges}).

In Step 3 we would hope to use the closeness of $H_n$ to the uniform distribution to convert the $dH_n(u)$ integrals to $du$ integrals. Then we write
$$
\frac{cd}{n} \int_0^1 \left\{F_c(u) - G_d(u)\right\}2\, du
$$
as a sum of three terms
$$
T_1 = \frac{d}{n} \int_0^1 c\left\{F_c(u)- u\right\}^2 \, du,
$$
$$
T_2 = \frac{c}{n} \int_0^1 d\left\{G_d(u)- u\right\}^2 \, du,
$$
and
$$
T_3 = -2\frac{\sqrt{cd}}{n}  \int_0^1 \sqrt{c} \left\{F_c(u)- u\right\} \cdot \sqrt{d}\left\{G_d(u)- u\right\}\, du.
$$
The integrals in $T_1$ and $T_2$ are both one sample Cram\'er-von Mises statistics so they are on the order 1. For any sequence $c=c_n$ such that $c_n/n\to 0$ the coefficient in front of $T_2$ is $o(1)$. So $T_2$ is negligible relative to $T_1$. The Cauchy-Schwarz inequality then shows $T_3$ is negligible relative to $T_2$.  There is a parallel argument when $c_n/m\to 1$.

Step 4 is not conjecture; its proof is straightforward from the assumptions of the Conjecture.  Steps 5 and 6 are exactly parallel to the arguments in \cite{lockhart1991overweight}. 

\bibliographystyle{plainnat}
\bibliography{Changepoints}
  \end{document}